\documentclass[a4paper,12pt,reqno]{amsart}
\usepackage{graphicx}
\usepackage{subcaption}
\usepackage{amsmath, amssymb, amsthm, epsfig}
\usepackage{latexsym}
\usepackage{url}
\usepackage{xcolor}
\usepackage{fullpage} 
\usepackage{setspace}
\usepackage{comment}

\usepackage{hyperref}

 \newcommand{\R}{\mathbb{R}}

\newcommand{\dr}{\text{\rm d}r}
\newcommand{\dx}{\text{\rm d}x}

\newcommand{\dd}{{\mathrm{d}}}

\newtheorem{theorem}{Theorem}

\newtheorem{lemma}[theorem]{Lemma}
\newtheorem{proposition}[theorem]{Proposition}
\newtheorem{corollary}{Corollary}

\theoremstyle{definition}

\title{Global and local maximizers for some Fourier extension estimates on the sphere}
\author{Valentina Ciccone and Mateus Sousa}
\date{\today}
\address[Valentina Ciccone]{Institute of Mathematics, Polish Academy of Sciences, Śniadeckich 8, 00-656 Warszawa, Poland}
\email{vciccone@impan.pl}
\address[Mateus Sousa]{BCAM - Basque Center for Applied Mathematics, Alameda de Mazarredo 14, 48009 Bilbao, Bizkaia,
Spain
}
\email{mcosta@bcamath.org}
\allowdisplaybreaks

\begin{document}

\begin{abstract}
In this note, we study maximizers for Fourier extension inequalities on the 
sphere. We prove that constant functions are local maximizers for the $L^p(\mathbb{S}^{d-1})$ to $L^p(\mathbb{R}^d)$ Fourier extension estimates  
in the same range of exponents $p$ for which they are global maximizers for the $L^2(\mathbb{S}^{d-1})$ to $L^p_{rad}L^2_{ang}(\mathbb{R}^d)$ mixed-norm Fourier extension inequalities. 
Moreover, in the case of low dimensions, we improve the range of exponents for which constant functions are known to be the unique global maximizers for the $L^2(\mathbb{S}^{d-1})$ to $L^p_{rad}L^2_{ang}(\mathbb{R}^d)$ mixed-norm Fourier extension estimate on the sphere, covering, for the case of dimensions $d=2,3$, the entire Stein--Tomas range.
This is achieved by establishing novel hierarchies between certain weighted norms of Bessel functions. 
\end{abstract}

\thanks{
V. Ciccone was partially supported by the Deutsche Forschungsgemeinschaft (DFG, German Research Foundation) under
Germany’s Excellence Strategy – EXC2047/1 390685813 as well as SFB 1060, and by the National Science Centre, Poland, grant Sonata Bis 2022/46/E/ST1/00036.
This project was partially carried out during V. Ciccone's research visits to BCAM - Basque Center for Applied Mathematics and to the University of the Basque Country UPV/EHU. Her visits have been partially supported by the Hausdorff Center for Mathematics in Bonn through the Global Math Exchange Program and BIGS,  and
by the grant PID2021-122156NB-I00 funded by MICIU/AEI/10.13039/501100011033 and FEDER, UE, and grant
IT1615-22 of the Basque Government.}

\thanks{ 
M. Sousa is is supported by grants RYC2022-038226-I and PID2020-113156GB-I00 funded by MICIU/AEI/10.13039/501100011033 and by ESF+, the Basque Government through the BERC 2022-2025 program, and through BCAM Severo Ochoa accreditation CEX2021-001142-S / MICIN / AEI / 10.13039/501100011033. }

\maketitle

\section{Introduction}

Let $d\geq 2$ be an integer, $J_\nu$ denote the Bessel function of the first kind of order $\nu$, and $k$ be a non-negative integer. 
It follows by the asymptotic behaviour of Bessel functions that the weighted $L^p$ norms 
$$\Lambda_{d,p}(k):=\bigg( \int_0^\infty |J_{\tfrac{d}{2}-1+k}(r)r^{-\tfrac{d}{2}+1}|^p r^{d-1} \dr \bigg)^{1/p} ~,$$
$$\Lambda_{d,\infty}(k):=\sup_{r\geq 0} |J_{\tfrac{d}{2}-1+k}(r)r^{ -\tfrac{d}{2}+1}|~,$$
are bounded whenever $\tfrac{2d}{d-1}<p\leq \infty$.

The problem of determining for which $k$ such weighted norms are maximized, which is a problem of independent interest in the theory of special functions, has been studied in \cite{COeSS19} in connection with certain mixed-norm sharp Fourier extension problems. In particular, in \cite{COeSS19} the authors have studied the problem of computing 
\begin{equation}\label{prob:Bessel} \tag{P1}
    \sup_{k\geq 0} \; \Lambda_{d,p}(k).
\end{equation}
The properties of Bessel functions (see e.g. \cite{St00}) guarantee that such supremum is a maximum. It has been shown in \cite{COeSS19} that such maximum is achieved at $k=0$ (and only at $k=0$) whenever $p$ is an even exponent and that the set of exponents for which the maximum is achieved at $k=0$ 
is open and it contains a neighborhood of infinity $(p_0(d),\infty ]$, providing some upper-bounds for $p_0(d)$. In particular, they obtained the following upper-bounds 
in low dimensions:
 $$p_0(2)\leq 6.76, \quad p_0(3)\leq 5.45, \quad p_0(4)\leq 5.53, \quad p_0(5)\leq 6.07, \quad p_0(6)\leq 6.82, $$
 $$p_0(7)\leq 7.70, \quad p_0(8)\leq 8.69, \quad p_0(9)\leq 9.78, \quad p_0(10)\leq 10.95,$$
and, more in general, they showed that
\begin{align}\label{COeSS_p0d_greater_dim}
p_0(d)\leq (\tfrac{1}{2}+o(1))d\log d.\end{align}
Problem \eqref{prob:Bessel} is related to several problems in sharp Fourier restriction theory. 

The Fourier restriction problem for the sphere asks for which pairs of exponents $(p,q)$ the inequality
\begin{align}\label{ineq:Fourier_extension_pq}
    \Vert \widehat{f\sigma} \Vert_{L^p(\mathbb{R}^d)} \leq C_{d,p,q} \Vert f\Vert_{L^q(\mathbb{S}^{d-1})}
\end{align}
holds. Here $\sigma=\sigma_{d-1}$ is the surface measure on $\mathbb{S}^{d-1}$ and $\widehat{f\sigma}$ is the Fourier transform of the measure $f\sigma$,
$$\widehat{f\sigma}(x)=\int_{\mathbb{S}^{d-1}}e^{-ix\cdot\xi}f(\xi)d\sigma(\xi).$$  The Fourier restriction problem has been fully solved only in dimension $d=2$ and for the case $q=2$ for which a complete answer is given by the Stein--Tomas inequality. A mixed-norm version of the problem has been studied in \cite{Ve92,Ve88} showing that the mixed-norm Fourier extension inequality 
\begin{equation}\label{vega}
    ||\widehat{f\sigma}||_{L^p_{rad}L^2_{ang}(\mathbb{R}^d)}\leq C_{d,p}||f||_{L^2(\mathbb{S}^{d-1})}
\end{equation}
holds when $\frac{2d}{d-1}<p$, where
\begin{align*}
    ||\widehat{f\sigma}||_{L^p_{rad}L^2_{ang}(\mathbb{R}^d)}=\left(\int_0^\infty\left(\int_{\mathbb{S}^{d-1}}|\widehat{f\sigma}(r\omega)|^2d\sigma(\omega)\right)^{p/2}r^{d-1}\dr\right)^{1/p}.
\end{align*}

The problem of determining the sharpest constant for \eqref{vega} has been studied in \cite{COeSS19}. Namely, in \cite{COeSS19} the authors have studied the problem of computing
\begin{equation}\label{prob:mixedNorm}\tag{P2}
\sup_{f\in L^2(\mathbb{S}^{d-1}),\; f\neq 0} \Phi^{\rm{mn}}_{p,d}(f),\quad \Phi^{\rm{mn}}_{p,d}(f):=\frac{||\widehat{f\sigma}||_{L^p_{rad}L^2_{ang}(\mathbb{R}^d)}}{||f||_{L^2(\mathbb{S}^{d-1})}}.
\end{equation}
It was observed in \cite{COeSS19} that the studying of such problem can be restricted to functions $f$ which are spherical harmonics. In other words
$$\sup_{f\in L^2(\mathbb{S}^{d-1}), \, f\neq 0}\Phi_{p,d}^{\rm{mn}}(f)=\sup_{Y_k,\, Y_k\neq 0}\Phi_{p,d}^{\rm{mn}}(Y_k),$$
 where $Y_k$ denotes a spherical harmonics of degree $k$.
Due to the identity
\begin{align}\label{FT_measure_on_sphere}
\widehat{Y_k\sigma}(x)=(2\pi)^{\frac{d}{2}}i^{-k}|x|^{-\frac{d}{2}+1}J_{\frac{d}{2}-1+k}(|x|)Y_k\big(\frac{x}{|x|}\big)
\end{align}
we have that
$$\Vert \widehat{Y_k\sigma}\Vert_{L^p_{rad}L^2_{ang}(\mathbb{R}^d)}=(2\pi)^{\tfrac{d}{2}}\bigg( \int_0^\infty |J_{\tfrac{d}{2}-1+k}(r)r^{-\tfrac{d}{2}+1}|^p r^{d-1} \dr\bigg)^{1/p}\Vert Y_k \Vert_{L^2(\mathbb{S}^{d-1})}.$$
Hence, the problem of establishing the sharpest constant for \eqref{vega}, namely \eqref{prob:mixedNorm}, is equivalent to the problem of determining for which non-negative integer $k$ the maximum in \eqref{prob:Bessel} is achieved.

Problem \eqref{prob:Bessel} has been addressed in \cite{COeSS19} by relating the integrals $\Lambda_{d,p}(k)$'s to integration on spheres using delta-calculus. Our approach, on the other hand, relies on some sharper estimates (with an improved constant) between weighted norms of Bessel functions 
inspired by those obtained in \cite[Lemma~2]{CG22} for the case of dimension $d=2$, see the forthcoming inequality \eqref{cg}.

Our first result lowers, for the case of low dimensions, the upper bounds for $p_0(d)$ established in \cite{COeSS19}, hence extending the ranges of exponents for which the maximum in \eqref{prob:Bessel} is achieved when $k=0$.
We use the notation $p_{\rm st}(d)$ to denote the Stein--Tomas endpoint exponent in dimension $d$, $p_{\rm st}(d):=\frac{2(d+1)}{(d-1)}$.
\begin{theorem}\label{wanted_thm}
 It holds that
     $$p_0(2)< 6, \quad  p_0(3)< 4, \quad p_0(4)< 3.48,  \quad p_0(5)< 3.50, $$ 
     $$ p_0(6)< 3.58, \quad p_0(7) < 3.7, \quad p_0(8) < 3.86, \quad p_0(9) < 4.06, \quad p_0(10) < 4.46. $$
     In particular, for $d=2,3$ this gives that $ p_0(d)< p_{\rm st}(d)$.
\end{theorem}
The fact that $ p_0(d)< p_{\rm st}(d)$ is of interest because constant functions are natural candidates to be extremizers for the full range of exponents of the Stein--Tomas Fourier extension inequality.  If this were true, then by H\"older's inequality, constant functions would be also maximizers for $\Phi_{p,d}^{\rm mn}$ when $p\geq p_{\rm{st}}(d)$. {This has been verified only when $p\geq 4$ is an even {integer} and $d\in\{3,4,5,6,7\}$ (see \cite{COeS15,Fo15,OSQ}), but it is open for all other cases.} In particular the case where $d=2$ has received a great deal of attention and many partial results have been achieved (see    \cite{CFOeST15,oliveira2019band,barker2020band,CG22,Lars23}), yet still remains unsolved. Hence, our result provides further evidence in this direction.

As mentioned above, Problems \eqref{prob:Bessel} and \eqref{prob:mixedNorm} are equivalent. Next, we observe that the same holds true if one considers the problem of finding extremizers among functions of the form $a Y_k\in L^2(\mathbb{S}^{d-1})$, with $a \in\mathbb{C}$ and $Y_k$ a spherical harmonic of degree $k$, for $L^p(\mathbb{S}^{d-1})$ to $L^p(\mathbb{R}^d)$ Fourier extension estimates. Namely, if one considers the problem of computing
\begin{equation}\label{prob:sphericHarmonics}\tag{P3}
\sup_{Y_k, \, Y_k\neq 0} \Phi_{p,d}(Y_k),\quad \Phi_{p,d}(f):=\frac{||\widehat{f\sigma}||_{L^p(\mathbb{R}^d)}}{||f||_{L^p(\mathbb{S}^{d-1})}}~.
\end{equation}
In fact,
$$\Vert \widehat{Y_k\sigma} \Vert_{L^p(\mathbb{R}^d)}=(2\pi)^{\tfrac{d}{2}}\bigg(\int_0^\infty |J_{\tfrac{d}{2}-1+k}(r)r^{-\tfrac{d}{2}+1}|^p r^{d-1}\dr \bigg)^{1/p}\Vert Y_k \Vert_{L^p(\mathbb{S}^{d-1})}~,$$
and, therefore,
$$\sup_{Y_k,\, Y_k\neq 0} \Phi_{p,d}(Y_k)=\sup_{f\in L^2(\mathbb{S}^{d-1}),\, f\neq 0} \Phi_{p,d}^{\rm{mn}}(f)=(2\pi)^{\tfrac{d}{2}}\sup_{k\geq 0}\Lambda_{d,p}(k).$$
In words, this simple observation asserts that the problem of computing the optimal constant for the mixed-norm Fourier extension inequality \eqref{vega} is equivalent to the problem of computing the optimal constant for the $L^p(\mathbb{S}^{d-1})$ to $L^p(\mathbb{R}^d)$ Fourier extension inequality when restricting to spherical harmonics. As an exemplifying application, by combining this observation with the fact, established in \cite{COeS15}, that constant functions are maximizers for the $L^p(\mathbb{S}^{d-1})$ to $L^p(\mathbb{R}^d)$ Fourier extension inequalities for $p\in 2\mathbb{Z}$, $p>\tfrac{2d}{d-1}$ one immediately gets that constant functions are also maximizers for $\Phi_{p,d}^{\mathrm{mn}}$ for all $p\in 2\mathbb{Z}$, $p>\tfrac{2d}{d-1}$. This was established in \cite{COeSS19} using delta-calculus. 

The following corollary is an immediate consequence of the above considerations.
\begin{corollary} 
    For all $p\in(p_0(d),\infty]$ 
    we have that
    \begin{align*}
        \sup_{f\in L^2(\mathbb{S}^{d-1}),\, f\neq 0} \Phi_{p,d}^{\rm{mn}}&(f)\leq \Phi_{p,d}^{\rm{mn}}(\mathbf{1})~,\\
        \sup_{Y_k,\, Y_k\neq 0}\Phi_{p,d}&(Y_k)\leq \Phi_{p,d}(Y_0)~.
    \end{align*}
That is, for all such $p$'s, constant functions are maximizers for \eqref{prob:mixedNorm} and \eqref{prob:sphericHarmonics}.
\end{corollary}
Note that the fact that constant functions are extremizers for \eqref{prob:sphericHarmonics} is a necessary condition for this to be the case also for the more general problem of computing
\begin{equation}\label{prob:LpLp}\tag{P4}
    \sup_{f\in L^p(\mathbb{S}^{d-1}),\, f\neq 0} \Phi_{p,d}(f).
\end{equation}
 
Extremizers for \eqref{prob:LpLp} are known only when $p$ is an even admissible exponent, in which case it has been shown in \cite{COeS15} that constant functions are maximizers, and when $p=\infty$ in which case the same conclusion holds \cite{FS22}. Except for these cases, even the question of the existence of global extremizers for Problem \eqref{prob:LpLp} is open, we refer to \cite{FS22} for recent results concerning existence of maximizers for Fourier extension inequalities on spheres. 
 Due to symmetry, constant functions would be natural candidate to be extremizers. Also, it was noted in \cite{CQ14} that constant functions are always solutions to the corresponding Euler--Lagrange equations for any admissible pair of exponents $(p,q)$ for the Fourier extension inequality \eqref{ineq:Fourier_extension_pq}, so, in particular, for any admissible pair $(p,p)$.
 
 A further intermediate step toward a solution of Problem \eqref{prob:LpLp} 
is to understand the behavior of local extremizers. Local extremizers have been studied before for the case of the endpoint Stein--Tomas inequalities in \cite{CFOeST15}, \cite{CS12}, and \cite{GN22} showing, respectively, that constant functions are local maximizers 
for such inequalities when $d=2$, when $d=3$, and when $2 \leq d \leq 60$.

Our second main result addresses this question by providing a further connection between Problems \eqref{prob:Bessel} and \eqref{prob:LpLp}. 

\begin{theorem}\label{thm:local}
    Let $d\geq 2$ and $p>\tfrac{2d}{d-1}$. Assume that the $L^p(\mathbb{S}^{d-1})$ to $L^p(\mathbb{R}^d)$ Fourier extension inequality holds and that the maximum in \eqref{prob:Bessel} is achieved at $k=0$. Then there exists $\delta >0$ such that whenever $\Vert f-\mathbf{1} \Vert_{L^p(\mathbb{S}^{d-1})}<\delta$, 
    \begin{equation}\label{inequa_local}
        \Phi_{p,d}(f)\leq \Phi_{p,d}(\mathbf{1}).
    \end{equation}
     That is, constant functions are local maximizers for \eqref{prob:LpLp}.
\end{theorem}

As an immediate consequence we have that constant functions are local maximizers for 
the $L^p(\mathbb{S}^{d-1})$ to $L^p(\mathbb{R}^d)$ Fourier extension inequality 
for all $p\in(p_0(d),\infty]$ for which the inequality holds and upper bounds on $p_0(d)$ are provided by Theorem \ref{wanted_thm} for the cases of dimensions $2\leq d \leq 10$, and, in general, by \eqref{COeSS_p0d_greater_dim} for greater dimensions. For example, this establishes that constant functions are local maximizers for such inequalities for all $p\geq 6$ in the case of dimension $d=2$, and  for all $p\geq 4$ in the case of dimension $d=3$.

The proof of Theorem \ref{thm:local} is contained in Section \ref{sec:proofLocal}, the proof of Theorem \ref{wanted_thm} is the content of Section \ref{sect:proof_thm_range}, while some auxiliary results about hierarchies between weighted norms of Bessel functions are presented in Section \ref{section:monoBessel}.

The topic of sharp spherical restriction has received much attention over the last years, in particular for the case of inequalities in the Stein--Tomas range \cite{FVV11,CS12,Fo15,COeS15,FLS16,Sh16,CFOeST15,oliveira2019band,barker2020band,OSQ,CG22,Lars23}. We refer to the survey \cite{NOeST23} for an up-to-date description of the state of the art.

\section{Hierarchies between weighted norms of Bessel functions}\label{section:monoBessel}

It is known that when $p\in 2\mathbb{N}$, $p>\frac{2d}{d-1}$, or when $p=\infty$ then
$$\frac{\Lambda_{d,p}(k)}{\Lambda_{d,p}(0)}< 1$$
for all positive integers $k$, see \cite{COeSS19}. In this section, we are interested in obtaining sharper estimates for such ratio, 
at least for certain values of the exponent $p$.

In this direction, for the case of dimension $d=2$ and exponent $p=6$ it has been shown in \cite{CG22} that
\begin{align}\label{cg}
\Lambda_{2,6}^6(k) < \frac{1}{3}\Lambda_{2,6}^6(0)
\end{align}
for all $k\geq 1$.

Moreover, in \cite{COeSS19} combining the identity
\begin{equation}\label{Linfty_0}
\Lambda_{d,\infty}(0)=\frac{1}{2^{\frac{d}{2}-1}\Gamma(\frac{d}{2})},\end{equation}
with a decreasing upper-bound (with respect to the order $k$) for $\Lambda_{d,p}(k)$, it has been shown that
$$\frac{\Lambda_{d,\infty}(k)}{\Lambda_{d,\infty}(0)}\leq \bigg( {\rm{L}}^6 \frac{2^{3d-6}\Gamma(\frac{d}{2})}{d^{3d-4}}\bigg)^{\frac{1}{3d+2}}$$
for all $k\geq 1$, where the constant ${\rm{L}}$ is defined as
\begin{equation}\label{Landau_L}
{\rm{L}}:=\sup_{\nu>0, \, r>0}|r^{1/3}J_\nu(r)|=0.785746...
\end{equation}
and it has been found by Landau \cite{La00}.

Our first result of this section establishes a hierarchy between the $\Lambda_{d,\infty}(k)$'s, hence determining the sharpest upper-bound on the ratio $\frac{\Lambda_{d,\infty}(k)}{\Lambda_{d,\infty}(0)}$.

\begin{proposition}\label{monotonicity_Linfty} For all positive integers $k$ it holds that
$$\Lambda_{d,\infty}(k-1)>\Lambda_{d,\infty}(k).$$
In particular,
$$\Lambda_{d,\infty}(k)\leq {\rm C}_{\infty} (d)\Lambda_{d,\infty}(0)$$
for all positive integers $k$, where ${\rm C}_\infty(d):=\frac{\Lambda_{d,\infty}(1)}{\Lambda_{d,\infty}(0)}$, and equality is attained if and only if $k=1$.   
\end{proposition}
\proof We begin with the case $d=2$. In such case
$\Lambda_{d,\infty}(k)=\sup_{r\geq 0}|J_k(r)|$. It has been shown in \cite{La00} that $\sup_{r> 0}|J_k(r)|$ is a strictly decreasing function of $k$. In particular, if we denote by $j'_{k,1}$ the first positive zero of $J_k'$ with $k$ a positive real number, then
\begin{align*}
    \sup_{r>0}|J_k(r)|=J_k(j_{k,1}'),
\end{align*}
and therefore
\begin{align*}
    \sup_{r>0}|J_k(r)|=J_k(j_{k,1}')> \sup_{r>0}|J_{k+1}(r)|=J_k(j_{k+1,1}').
\end{align*}
As $\sup_{r\geq 0}|J_0(k)|=J_0(0)=1>J_1(j_{1,1}')$ the claim in the statement is verified for the case $d=2$.

We turn to the case of $d\geq 3$. In these cases $\Lambda_{d,\infty}(k)=\sup_{r\geq 0}|r^{-\tfrac{d}{2}+1}J_{\tfrac{d}{2}-1+k}|$. We start by observing that as $r^{-\tfrac{d}{2}+1}$ is a strictly decreasing function of $r$ and $\sup_{r>0}|J_\nu(r)|=J_\nu(j_{\nu,1}')$ it holds that 
\begin{align*}
    \sup_{r>0}|r^{-\tfrac{d}{2}+1}J_{\tfrac{d}{2}-1+k}(r)|=\sup_{0 < r < j_{d/2-1+k,1}'} |r^{-\tfrac{d}{2}+1}J_{\tfrac{d}{2}-1+k}(r)|.
\end{align*}
Hence to conclude it would be enough to show that
\begin{align}\label{proof_monotonicity_infty_what_wanted}
    J_\nu(r)>J_{\nu+1}(r) \qquad \text{for all } r\in(0,j_{\nu,1}').
\end{align}
We recall that $j_{\nu,1}'< j_{\nu+1,1}'$, see e.g. \cite{Wa66}. In particular, $J_\nu, \; J_{\nu+1}, \, J_\nu',\, J_{\nu+1}'$ are strictly positive in $(0,j'_{\nu,1})$, and $J_{\nu}(j_{\nu,1}')>J_{\nu+1}(j_{\nu,1}')$. Hence to prove \eqref{proof_monotonicity_infty_what_wanted} it suffices to show that there exist no $\overline{r}\in (0,j_{\nu,1}')$ such that $J_\nu(\overline{r})=J_{\nu+1}(\overline{r})$. We argue by contradiction. Consider the recursive relations for Bessel functions
\begin{align}\label{first_recursive}
    \tfrac{2\nu}{r}J_\nu(r)=J_{\nu-1}(r)+J_{\nu +1}(r),
\end{align}
\begin{align}\label{second_recursive}
    2J'_\nu(r)= J_{\nu-1}(r)-J_{\nu+1}(r).
\end{align}
By taking the sum of \eqref{first_recursive} and \eqref{second_recursive} we obtain the identity
$$J_\nu'(r)=J_{\nu-1}(r)-\tfrac{(\nu)}{r}J_{\nu}(r),$$
and shifting $\nu \mapsto \nu +1$ we obtain
$$J_{\nu+1}'(r)=J_{\nu}(r)-\tfrac{(\nu+1)}{r}J_{\nu+1}(r).$$
Assume there exist $\overline{r}\in (0,j_{\nu,1}')$ such that $J_\nu(\overline{r})=J_{\nu+1}(\overline{r})$. Evaluating the last display at $\overline{r}$ we get
$$J_{\nu+1}'(\overline{r})=(1-\tfrac{(\nu+1)}{\overline{r}})J_{\nu+1}(\overline{r})$$
and, as $J_{\nu+1}'$ and $J_{\nu+1}$ are strictly positive on $(0,j_{\nu,1}')$, we have that necessarily
\begin{align}\label{proof_monotonicity_infty_first_condition_contradiction}
    \overline{r}>\nu+1.
\end{align}
Next, we take the difference between \eqref{second_recursive} and \eqref{first_recursive} obtaining the identity
$$J_{\nu}(r)'=\tfrac{\nu}{r}J_\nu(r)-J_{\nu+1}(r).$$
Evaluating it at $\overline{r}$ we get that
$$J_\nu'(\overline{r})=(\tfrac{\nu}{\overline{r}}-1)J_\nu(\overline{r})$$
and, as both $J_\nu'$ and $J_\nu$ are strictly positive on $(0,j_{\nu,1}')$, we have that necessarily
$$\nu> \overline{r}.$$
Comparing this with \eqref{proof_monotonicity_infty_first_condition_contradiction} yields the contradiction.
\qed

The values of $\Lambda_{d,\infty}(1)$ can be computed using Mathematica. For the case of $2\leq d\leq 10$ one obtains, with 6 significant figures (s.f.),
\begin{align}\begin{split}\label{Lambda_infty_one}
\Lambda_{2,\infty}(1)=0.581865, \quad \Lambda_{3,\infty}(1)&=0.348023, \quad\Lambda_{4,\infty}(1)=0.179963, \\
\Lambda_{5,\infty}(1)=0.0830013, \quad \Lambda_{6,\infty}(1)&=0.0348492, \quad  \Lambda_{7,\infty}(1)=0.0135129, \\
\Lambda_{8,\infty}(1)=0.00489072, \quad \Lambda_{9,\infty}(1)&=0.00166575, \quad \Lambda_{10,\infty}(1)=0.000537364 .
\end{split}
\end{align}
By combining them with \eqref{Linfty_0} one can obtain a numerical evaluation for ${\rm C}_\infty(d)$.

Our second observation is for the case of exponent $p=4$ and dimensions $3\leq d \leq 10$. 

\begin{proposition}\label{proposition_p4}
    Let $3\leq d \leq 10$. Then
$$ \Lambda_{d,4}(k)\leq {\rm C}_4(d)\, \Lambda_{d,4}(0)$$
holds for all positive integers $k$, where ${\rm C}_4(d):=\tfrac{\Lambda_{d,4}(1)}{\Lambda_{d,4}(0)}  \, <1.$ Equality is attained if and only if $k=1$.
\end{proposition}

To prove Proposition \ref{proposition_p4}, we rely on the following upper-bound; see also \cite{GN22}.

\begin{lemma}\label{lemma_upperbound}
    Let $d\geq 2$. In the range of exponents $\tfrac{6d-2}{3d-4}<p<\tfrac{12d+4}{3d-4}$ it holds that
 $$\Lambda_{d,p}^p(k)\leq {\rm L}^{p-2}\frac{\Gamma(\lambda)\Gamma(\tfrac{d}{2}-1+k+\frac{1-\lambda}{2})}{2^\lambda \Gamma(\tfrac{1+\lambda}{2})^2\Gamma(\tfrac{d}{2}-1+k+\tfrac{1+\lambda}{2})}$$
    for all positive integers $k$, where $\lambda=p(\tfrac{d}{2}-\tfrac{2}{3})-d+\tfrac{1}{3}$.
\end{lemma}
To establish the upper-bound in the lemma, we rely on the following identity, which can be found in \cite[Equation~6.574-2]{GR14}
    \begin{align}\label{gradshteyntable}
        \int_0^\infty J_\nu^2(r)r^{-\lambda}\dr= \frac{\Gamma(\lambda)\Gamma(\nu+\tfrac{(1-\lambda)}{2})}{2^\lambda\Gamma\big( \tfrac{1+\lambda}{2} \big)^2 \Gamma\big( \nu + \tfrac{1+\lambda}{2}\big)}
    \end{align}
for $0<\lambda<2\nu+1$.
\proof
We use \eqref{Landau_L} to obtain the upper-bound
$$\Lambda_{d,p}^p(k)\leq {\rm L}^{p-2}\int_0^\infty J_{\tfrac{d}{2}-1+k}^2(r)r^{-p\big(\tfrac{d}{2}-\tfrac{2}{3}\big)+d-\tfrac{1}{3}}\dr.$$
By applying identity \eqref{gradshteyntable} to the right hand side of the last display we further obtain
$$\Lambda_{d,p}^p(k)\leq {\rm L}^{p-2}\frac{\Gamma(\lambda)\Gamma(\tfrac{d}{2}-1+k+\frac{1-\lambda}{2})}{2^\lambda \Gamma(\tfrac{1+\lambda}{2})^2\Gamma(\tfrac{d}{2}-1+k+\tfrac{1+\lambda}{2})}$$
where $\lambda=p(\tfrac{d}{2}-\tfrac{2}{3})-d+\tfrac{1}{3}$. Such upper-bound holds whenever $0<\lambda< 2\big(\tfrac{d}{2}-1 +k\big) +1.$ In particular, for a fixed dimension $d\geq 2$ the upper-bound holds for all positive integers $k$ whenever $\tfrac{6d-2}{3d-4}<p<\tfrac{12d+4}{3d-4}$.
\qed

Note that both the case of $p=4$ and the case of Stein--Tomas endpoint 
$p_{\rm st}(d)$ are included in the range of exponents covered by Lemma \ref{lemma_upperbound}. Also, note that, for a fixed exponent $p$ and a fixed dimension $d$, the above upper bound is a decreasing function of $k$. Throughout, we use the notation $U_{d,p}(k)$ to denote the upper bound for $\Lambda_{d,p}^p(k)$ in Lemma \ref{lemma_upperbound}.

\textit{Proof of Proposition \ref{proposition_p4}.} \,
We compare the upper bound $U_{d,4}(k)$ for $\Lambda_{d,4}^4(k)$ established in Lemma \ref{lemma_upperbound} with a (lower) estimate for $\Lambda_{d,4}^4(1)$. To this end, we rely on Mathematica to evaluate the integrals
\begin{align}\label{Lambda_4_1_truncated}
\int_0^{40} |J_{\frac{d}{2}}(r)r^{-\tfrac{d}{2}+1}|^4 r^{d-1}\dr
\end{align}
for $3\leq d\leq 10$ obtaining, respectively, the following values (with 6 s.f.)
\begin{align*}
0.144681 \quad 0.0337263 & \quad 0.00661348 \quad 0.00107217 \\ 0.000146318 \quad 0.0000171549 &\quad 1.75867 \times 10^{-6} \quad 1.59953 \times 10^{-7}.
\end{align*}
By comparison, one can see that $U_{d,4}(k)< \Lambda_{d,4}^4(1)$ for all integers $k\geq 2$ when $d\in\lbrace 5,6,7,8 
\rbrace$, for all integers $k\geq 3$ when $d\in \lbrace 4,9,10 \rbrace$, and for all integers $k\geq 5$ when $d=3$. We check the remaining cases separately. We rely on Mathematica to evaluate the integrals
$$\int_0^{200} |J_{\tfrac{d}{2}-1+k}(r) r^{-\tfrac{d}{2}+1}|^4 r^{d-1}\dr$$
for the cases of interest obtaining, for the case of $d=3$ and $k=2,3,4$, the values (6 s.f.)
$$ 0.0992828 \quad 0.0757045 \quad 0.0615859, $$
respectively, and for the cases 
$k=2$ and $d=4,9,10$, the values (6 s.f.) 
$$0.0172602 \quad  4.70782\times 10^{-7} \quad 4.00184\times 10^{-8},$$
respectively. Then, we use the estimate
\begin{align}\label{upper_bound_power_half}
|J_\nu(r)| \leq r^{-1/2}
\end{align}
which holds for all $\nu\geq \tfrac{1}{2}$ and $r\geq \tfrac{3}{2}\nu$ (see \cite[Lemma~8]{COeSS19} and \cite[Theorem~3]{Kr14}) to upper bound the tails 
obtaining that
$$\int_{200}^\infty |J_{\tfrac{d}{2}-1+k}(r) r^{-\tfrac{d}{2}+1}|^4 r^{d-1}\dr\leq \frac{200^{-d+2}}{d-2}.$$
Hence, by comparison, we see that also for these cases it holds that $\Lambda_{d,4}^4(k) < \Lambda_{d,4}^4(1)$.  As it is known from \cite{COeS15,COeSS19} that $\Lambda_{d,4}^4(1)< \Lambda_{d,4}^4(0)$ when $d\geq 3$, the result in the statement follows.
\qed

To evaluate ${\rm C}_4(d)$ one can rely on the identity
$$\int_0^\infty |J_\nu(r)|^4 r^{-2\nu +1} \dr =\frac{\Gamma(\nu) \Gamma(2\nu)}{2\pi \Gamma(\nu + \tfrac{1}{2})^2 \Gamma(3\nu)}, $$
which can be found, for example, in \cite[Lemma~7]{COeSS19} (see also \cite[Equation~6.5793-3]{GR14}) and which provides an explicit expression for $\Lambda_{d,4}^4(0)$, together with 
 a numerical estimates for $\Lambda_{d,4}^4(1)$. 

Our last result of this section is for the case of the Stein--Tomas endpoint, 
$p_{\rm st}=p_{\rm{st}}(d)$. 
\begin{proposition}\label{proposition_p_st}
    Let $4\leq d \leq 10$. Then the following inequality holds for all positive integers $k$
    $$\Lambda_{d,p_{\rm{st}}}(k) \,  \leq \, {\rm C}_{p_{\rm st}}(d)\; \Lambda_{d,p_{\rm{st}}}(0),$$
   where ${\rm C}_{p_{\rm st}}(d):=\tfrac{\Lambda_{d,p_{\rm st}}(1)}{\Lambda_{d,p_{\rm st}}(0)}  \, <1$. Equality is attained if and only if $k=1$.
\end{proposition}

\proof  
We compare the upper bound $U_{d,p_{\rm st}}(k)$ for $\Lambda_{d,p_{\rm st}}^{p_{\rm st}}(k)$ established in Lemma \ref{lemma_upperbound} with a (lower) estimate for $\Lambda_{d,p_{\rm st}}^{p_{\rm st}}(1)$. To this end, we rely on Mathematica to evaluate the integrals
\begin{align}\label{truncated_pst_k1}
\int_0^{50} |J_{\frac{d}{2}}(r)r^{-\tfrac{d}{2}+1}|^{p_{\rm st}} r^{d-1}\dr
\end{align}
for $4\leq d\leq 10$ obtaining, respectively, the following values (with 6 s.f.)
$$0.143391 \quad 0.131693 \quad 0.118941 \quad 0.10719 \quad 0.0969753 \quad 0.088279 \quad 0.0807943. $$
By comparison, one can see that $U_{d,p_{\rm st}}(k)< \Lambda_{d,p_{\rm st}}^{p_{\rm st}}(1)$ for all integers $k\geq 3$ when $d\in\lbrace 5,6,7,8,9,10 \rbrace$, and for all integers $k\geq 4$ when $d=4$.
We check the remaining cases separately. We use Mathematica to evaluate the integral
$$\int_0^{200} |J_{\tfrac{d}{2}-1+k}(r) r^{-\tfrac{d}{2}+1}|^{p_{\rm st}} r^{d-1}\dr$$
obtaining for the cases $k=2$ and $5\leq d \leq 10$ the values (6 s.f.)
$$0.0998066 \quad 0.0938562 \quad 0.0875322\quad 0.0814907 \quad 0.075952 \quad 0.0709569 $$
and for the cases $d=4$ and $k=2,3$ the values (6 s.f.)
$$0.103492 \quad 0.080522.$$
We use the estimate \eqref{upper_bound_power_half} to upper bound the tail of the integrals obtaining
$$\int_{200}^\infty |J_{\tfrac{d}{2}-1+k}(r) r^{-\tfrac{d}{2}+1}|^{p_{\rm st}} r^{d-1}\dr \leq   \frac{1}{200}.$$
Hence, by comparison, it follows that $\Lambda_{4,p_{\rm st}}^{p_{\rm st}}(2)< \Lambda_{4,p_{\rm st}}^{p_{\rm st}}(1)$.

We are left to show that $\Lambda_{d,p_{\rm st}}^{p_{\rm st}}(1)< \Lambda_{d,p_{\rm st}}^{p_{\rm st}}(0)$ whenever $4\leq d \leq 10$. The cases of $d=4,5$ have already be verified in \cite{COeSS19}. To verify the remaining cases $6\leq d \leq 10$ we compare the  bound for $\Lambda_{d,p_{\rm st}}^{p_{\rm st}}(1)$  obtained by combining the numerical evaluation of the truncated integral \eqref{truncated_pst_k1} and an upper bound for the tail obtained using \eqref{upper_bound_power_half}  with
a (lower) estimate for $\Lambda_{d,p_{\rm st}}^{p_{\rm st}}(0)$. To this end we numerically evaluate the integral
$$\int_0^{50} |J_{\frac{d}{2}+1}(r)r^{-\tfrac{d}{2}+1}|^{p_{\rm st}} r^{d-1}\dr$$ for $6\leq d \leq 10$ obtaining the values (6 s.f.)
$$0.173201 \quad 0.147926 \quad 0.1286 \quad 0.113331 \quad 0.101086.$$
By comparison, we see that $\Lambda_{d,p_{\rm st}}^{p_{\rm st}}(1)< \Lambda_{d,p_{\rm st}}^{p_{\rm st}}(0)$ also for $6\leq d \leq 10$ hence concluding the proof.
\qed

\section{Proof of Theorem \ref{wanted_thm}}\label{sect:proof_thm_range}
\subsection{Case $d=2$}
We combine the 
estimate \eqref{cg} from \cite{CG22} and the 
estimate in Proposition \ref{monotonicity_Linfty} (here for the case $d=2$) with the interpolation strategy utilized in \cite{COeSS19}. Let $p\geq 6$ and $k$ be a positive integer. It follows from H{\" o}lder's inequality that
$$\Lambda_{2,p}(k)\leq \Lambda_{2,6}(k)^{6/p}\Lambda_{2,\infty}(k)^{1-6/p}.$$
Using \eqref{cg} and the sharp estimate from Proposition \ref{monotonicity_Linfty} we further obtain that
$$\Lambda_{2,p}(k)\leq \frac{1}{3^{1/p}} \Lambda_{2,6}(0)^{6/p}\Lambda_{2,\infty}(1)^{1-6/p}. $$
We need the following lower bound on $\Lambda_{d,p}(0)$ which has been established in \cite[Equation~4.8]{COeSS19}
\begin{align}\label{lower_bound_Lambdak0}
    \Lambda_{d,p}(0)> \frac{(2^{d-1}(\tfrac{d}{2})^{d/2})^{1/p}}{2^{d/2-1}\Gamma(\tfrac{d}{2})} \bigg( \frac{\Gamma(p+1)\Gamma(\tfrac{d}{2})}{\Gamma(p+\tfrac{d}{2}+1)} \bigg)^{1/p}.
\end{align}
Then, we rely on standard numerical evaluation to determine for which $p\geq 6$ it holds that 
$$ \frac{1}{3^{1/p}} \Lambda_{2,6}(0)^{6/p}\Lambda_{2,\infty}(1)^{1-6/p} \leq 
2^{1/p}\frac{\Gamma(p+1)}{\Gamma(p+2)}.$$
We obtain that such inequality is satisfied for all $p\geq 6$. Hence, $p_0(2) < 6$ as claimed.

\subsection{Case $d\geq 3$}
We proceed in two steps.
First, we combine the 
estimates in Proposition \ref{proposition_p4} and Proposition \ref{monotonicity_Linfty} with the interpolation strategy utilized in \cite{COeSS19}. This will establish the upper bound on $p_0(d)$ in the statement of Theorem \ref{wanted_thm}
 for the cases of $d=3,9,10$. Second, we use the 
 estimates in Proposition \ref{proposition_p4} and Proposition \ref{proposition_p_st} and interpolation to establish the upper bound on $p_0(d)$ 
 for the cases of $d=4,5,6,7,8$. 

 \subsubsection*{Step 1} Let $p\geq 4$ and $k$ be a positive integer. It follows from H{\" o}lder's inequality that
 $$\Lambda_{d,p}(k)\leq \Lambda_{d,4}(k)^{4/p}\Lambda_{d,\infty}(k)^{1-4/p}. $$
 Using the sharp estimate from Proposition \ref{monotonicity_Linfty} and Proposition \ref{proposition_p4} we further obtain that
 $$\Lambda_{d,p}(k)\leq \Lambda_{d,4}(1)^{4/p}\Lambda_{d,\infty}(1)^{1-4/p}. $$
 Then, we compare the right-hand side of the last display with the lower bound for $\Lambda_{d,p}(0)$ in equation \eqref{lower_bound_Lambdak0} to determine for which $p\geq 4$ the following inequality is satisfied
$$ \Lambda_{d,4}(1)^{4/p}\Lambda_{d,\infty}(1)^{1-4/p} \leq \frac{(2^{d-1}(\tfrac{d}{2})^{d/2})^{1/p}}{2^{d/2-1}\Gamma(\tfrac{d}{2})} \bigg( \frac{\Gamma(p+1)\Gamma(\tfrac{d}{2})}{\Gamma(p+\tfrac{d}{2}+1)} \bigg)^{1/p}.$$
 We use the numerical values for $\Lambda_{d,\infty}(1)$ in \eqref{Lambda_infty_one} and the bound for $\Lambda_{d,4}(1)$ obtained by combining the numerical evaluation for the truncated integral in \eqref{Lambda_4_1_truncated} with an upper bound for the tail obtained using \eqref{upper_bound_power_half}. Via a standard numerical evaluation, we obtain that $\Lambda_{d,p}(k) < \Lambda_{d,p}(0)$ for all $p \geq 4$ for the cases of dimensions $d=3,4,5,6,7,8$, for all $p \geq 4.06$ for the case of $d=9$, and for all $p\geq 4.46$ for the case of $d=10$.

\subsubsection*{Step 2} Let $4\leq d \leq 8$,  $p_{\rm st}(d)\leq p \leq 4$ and $k$ be a positive integer. It follows from H{\" o}lder's inequality that 
$$\Lambda_{d,p}(k)\leq \Lambda_{d,p_{\rm st}}(k)^{(1-\theta)}\Lambda_{d,4}(k)^{\theta}, $$
with $\theta:=\tfrac{4}{p}\tfrac{(p-p_{\rm st})}{(4- p_{\rm st})}$.
Using the estimates of Proposition \ref{proposition_p4} and Proposition \ref{proposition_p_st} we further obtain that
$$\Lambda_{d,p}(k)\leq \Lambda_{d,p_{\rm st}}(1)^{(1-\theta)}\Lambda_{d,4}(1)^{\theta}.$$
As before, we bound $\Lambda_{d,4}(1)$ by combining the numerical evaluation for the truncated integral in \eqref{Lambda_4_1_truncated} with an upper bound for the tail obtained using \eqref{upper_bound_power_half} and we proceed analogously for $\Lambda_{d,p_{\rm st}}(1)$.
Then, we compare this upper bound with the lower bound for  $\Lambda_{d,p}(0)$ in equation \eqref{lower_bound_Lambdak0} to determine, for a fixed $4\leq d \leq 8$ ,  for which $p_{\rm st}(d)\leq p \leq 4$ the former is greater than the latter. 
We obtain that $\Lambda_{d,p}(k) < \Lambda_{d,p}(0)$ for all 
$p\geq p(d)$ with
$$p(4)=3.48 ,\quad p(5)=3.5 ,\quad p(6)= 3.58 \quad, p(7)= 3.7, \quad p(8)= 3.86  . $$
\qed

\section{Proof of Theorem \ref{thm:local}}\label{sec:proofLocal}
Consider the deficit functional
\begin{align*}
\zeta_{p}[f]= \Phi_{p,d}(\mathbf{1})^p\|f\|_{L^p(\mathbb{S}^{d-1})}^p-\|\widehat{f\sigma}\|_{L^p(\mathbb{R}^d)}^p.
\end{align*}
Inequality \eqref{inequa_local}  is equivalent to 
\begin{equation}\label{eq:local_deficit}
    \zeta_p[f]\geq 0,
\end{equation}
therefore it is enough to prove that there is a $\delta>0$ such that $\zeta_p[f]>0$ when $\|f-\mathbf{1}\|_{L^p(\mathbb{S}^{d-1})}<\delta$ and $f$ is not constant, which we proceed to do. 

We recall that here $p>2$ and we are assuming that the Fourier extension operator is bounded from $L^p(\mathbb{S}^{d-1})$ to $L^p(\mathbb{R}^d)$. 
We compute
\begin{align}\label{eq:expansion1}\begin{split}
\int_{\mathbb{R}^d}|\widehat{\mathbf{1}\sigma}(x)+\varepsilon\widehat{g\sigma}(x)|^p \dx =&\int_{\R^d}|\widehat{\mathbf{1}\sigma}(x)|^p \dx+  p\varepsilon\int_{\R^d}|\widehat{\mathbf{1}\sigma}(x)|^{p-2}\mathfrak{R}(\widehat{\mathbf{1}\sigma}(x){\widehat{g\sigma}(x)})\dx \\
    &+\frac{p(p-2)\varepsilon^2}{4} \mathfrak{R} \int_{\R^d}|\widehat{\mathbf{1}\sigma}(x)|^{p-4}(\widehat{\mathbf{1}\sigma}(x){\widehat{g\sigma}(x)})^2 \dx\\
    &+\frac{p\varepsilon^2}{4}\int_{\R^d}|\widehat{\mathbf{1}\sigma}(x)|^{p-2}|\widehat{g\sigma}(x)|^2 \dx \\ &+
    o(\varepsilon^2).
    \end{split}
\end{align}
and 
\begin{align}\label{eq:expansion2}\begin{split}
    \int_{\mathbb{S}^{d-1}}|{\mathbf{1}}+\varepsilon{g(x)}|^p\dd\sigma(x) =&\|{\mathbf{1}}\|_{L^p(\mathbb{S}^{d-1})}^p+  p\varepsilon\int_{\mathbb{S}^{d-1}}\mathfrak{R}(g(x))\dd\sigma(x) \\
    & +\frac{p(p-2)\varepsilon^2}{4}\mathfrak{R}\int_{\mathbb{S}^{d-1}}{g(x)}^2 \dd\sigma(x) \\
    &+\frac{p\varepsilon^2}{4}\int_{\mathbb{S}^{d-1}}|{g(x)}|^2 \dd\sigma(x) + 
    o(\varepsilon^2).
    \end{split}
\end{align}
We take $f$ to be of the form $f=\mathbf{1}+\varepsilon g$, with $0<\varepsilon \leq \delta$ and $\Vert g \Vert_{L^p(\mathbb{R}^d)}=1$.
By applying \eqref{eq:expansion1} and \eqref{eq:expansion2} 
one has
\begin{align}\label{eq:deficit}
    \begin{split}
        \zeta_p[f]=& \zeta_p[\mathbf{1}+\varepsilon g]  
        =\,p\varepsilon\, \left(\Phi_{p,d}(\mathbf{1})^p \int_{\mathbb{S}^{d-1}} \mathfrak{R} (g(x)) \dd\sigma(x)  -\int_{\R^d}| \widehat{\mathbf{1}\sigma}(x)|^{p-2}\mathfrak{R} (\widehat{\mathbf{1}\sigma}(x) {\widehat{g\sigma}(x)}) \dx\right) \\
        &+\frac{\varepsilon^2}{4}p(p-2)\left(\Phi_{p,d}(\mathbf{1})^p\mathfrak{R} \int_{\mathbb{S}^{d-1}}{g(x)}^2 \dd\sigma(x) - \mathfrak{R}\int_{\R^d}|\widehat{\mathbf{1}\sigma}(x)|^{p-2}(\widehat{g\sigma}(x))^2 \dx\right)\\
        &+\frac{\varepsilon^2}{4}p\left(\Phi_{p,d}(\mathbf{1})^p\int_{\mathbb{S}^{d-1}}|{g(x)}|^2 \dd\sigma(x)-\int_{\R^d}|\widehat{\mathbf{1}\sigma}(x)|^{p-2}|\widehat{g\sigma}(x)|^2 \dx\right)+ 
        o(\varepsilon^2).
    \end{split}
\end{align}
Furthermore, due to the aforementioned observation that $\mathbf{1}$ is a critical point of $\Phi_{p,d}$, the first order terms in $\varepsilon$ of \eqref{eq:deficit} all vanish. To deal with the second order terms, we use the fact that $L^p(\mathbb{S}^{d-1})\subset L^2(\mathbb{S}^{d-1})$ since $p>2$ in order to expand $g$ in spherical harmonics. For that purpose we choose for each $k$ an orthonormal basis $\{Y_{j,k}\}_{j}$ of $\mathcal{H}_k^d$ where each $Y_{j,k}$ is a real-valued spherical harmonic of degree $k$. Then
\begin{align*}\begin{split}
        &g=\sum_{k,j} a_{j,k}Y_{j,k}. 
    \end{split}
\end{align*}
By combining identity \eqref{FT_measure_on_sphere} with the observation that the first order terms vanish at \eqref{eq:deficit} we can integrate in polar coordinates to obtain
\begin{align}\label{eq:deficit2}
    \begin{split}
        \zeta_p[f]
        =&\frac{\varepsilon^2}{4}p(p-2)\bigg(\Phi_{p,d}(\mathbf{1})^p\sum_{k,j}\mathfrak{R}(a_{j,k})^2 \\
        & -(2\pi)^{pd/2}\sum_{k,j} (-1)^k\mathfrak{R}(a_{j,k})^2\int_0^\infty|J_{\frac{d}{2}-1}(r)|^{p-2}|J_{\frac{d}{2}-1+k}(r)|^2r^{d-1-p(1-\frac{d}{2})} \dd{r}\bigg)\\
        &+\frac{\varepsilon^2}{4}p\sum_{k,j}|a_{j,k}|^2\bigg(\Phi_{p,d}(\mathbf{1})^p -(2\pi)^{pd/2}\int_0^\infty|J_{\frac{d}{2}-1}(r)|^{p-2}|J_{\frac{d}{2}-1+k}(r)|^2r^{d-1-p(1-\frac{d}{2})} \dd{r}\bigg) \\
        & + o(\varepsilon^2).   
    \end{split}
\end{align}
Lastly, using  H\"older's inequality and the fact that by hypothesis $\Lambda_{d,p}(k)<\Lambda_{d,p}(0)$ for all positive integers $k$ we observe that
\begin{align*}
    \int_0^\infty |& J_{\frac{d}{2}-1}(r)|^{p-2} |J_{\frac{d}{2}-1+k}(r)|^2 r^{d-1-p(1-\frac{d}{2})} \dd{r} \\
    & < \bigg(\int_0^\infty |J_{\tfrac{d}{2}-1}(r)|^{p} r^{d-1-p(1-\frac{d}{2})} \dd{r} \bigg)^{(p-2)/p} \bigg( \int_0^\infty |J_{\tfrac{d}{2}-1+k}(r)|^{p} r^{d-1-p(1-\frac{d}{2})} \dd{r}  \bigg)^{2/p}\\
    & < \int_0^\infty |J_{\tfrac{d}{2}-1}(r)|^{p} r^{d-1-p(1-\frac{d}{2})} \dd{r} = (2\pi)^{-pd/2}\Phi_{p,d}(\mathbf{1})^p,
\end{align*}
hence concluding the proof of Theorem \ref{thm:local}.

 \section*{Acknowledgements}
V. Ciccone thanks BCAM - Basque Center for Applied Mathematics, the University of the Basque Country UPV/EHU, and her hosts, Odysseas Bakas and Ioannis Parissis, for their warm hospitality during her research visits in Bilbao. 

The authors thank Christoph Thiele and Emanuel Carneiro for valuable comments during the preparation of this manuscript.


\begin{thebibliography}{100}
 
 




\bibitem[BTZK20]{barker2020band}
\textsc{J. Barker, C. Thiele, and P. Zorin-Kranich},
\newblock {\it Band-limited maximizers for a {F}ourier extension inequality on the circle {I}{I}},
Exp. Math. {\bf 32} (2023), no. 2, 280–293.

\bibitem[Be23]{Lars23}
\textsc{L. Becker},
\newblock{\it Sharp Fourier extension for functions with localized support on the circle},
\newblock {Rev. Mat. Iberoam. (2024).} 


\bibitem[CFOeST15]{CFOeST15}
\textsc{E. Carneiro, D. Foschi, D. Oliveira e Silva and C. Thiele},
\newblock {\it A sharp trilinear inequality related to Fourier restriction on the circle},
\newblock Rev. Mat. Iberoam. {\bf 33} (2017), no. 4, 1463--1486.

\bibitem[{COeS15}]{COeS15}
\textsc{E. Carneiro and D. Oliveira e Silva},
\newblock {\it Some sharp restriction inequalities on the sphere,}
\newblock Int. Math. Res. Not. IMRN (2015), no.~17, 8233--8267.

\bibitem[{COeSS19}]{COeSS19}
\textsc{E. Carneiro, D. Oliveira e Silva and M. Sousa},
\newblock {\it Sharp mixed norm spherical restriction},
\newblock  Adv. Math. {\bf 341} (2019), 583--608. 

\bibitem[CG24]{CG22}
\textsc{V. Ciccone and F. Gon\c calves}, 
\newblock{\it Sharp Fourier extension on the circle under arithmetic constraints}. J. Funct. Anal. \textbf{286} (2024), no.2, Paper No. 110219.


\bibitem[{CS12} ]{CS12} 
\textsc{M. Christ and S. Shao},
\newblock {\it Existence of extremals for a Fourier restriction inequality},
\newblock Anal. PDE. {\bf 5} (2012), no.~2, 261--312.

\bibitem[{CQ14} ]{CQ14} 
\textsc{M. Christ and R. Quilodr\'{a}n},
\newblock {\it  Gaussians rarely extremize adjoint Fourier restriction inequalities for paraboloids}. Proc. Amer. Math. Soc. {\bf 142} (2014), no. 3, 887--89

\bibitem[{FVV11}]{FVV11}
\textsc{L. Fanelli, L. Vega and N. Visciglia}, 
\newblock {\it On the existence of maximizers for a family of restriction theorems},
\newblock Bull. Lond. Math. Soc. {\bf 43} (2011), no.~4, 811--817.


\bibitem[FS24]{FS22}
\textsc{T. Flock and B. Stovall},
\newblock {\it  On Extremizing Sequences for Adjoint Fourier Restriction},
\newblock Adv. Math. \textbf{453} (2024), Paper No. 109854.

\bibitem[FLS16]{FLS16}
\textsc{R. Frank, E. Lieb and J. Sabin},
\newblock {\it Maximizers for the Stein--Tomas inequality}, 
\newblock Geom. Funct. Anal. {\bf 26} (2016), no. 4, 1095-1134.

\bibitem[{Fo15}]{Fo15}
\textsc{D. Foschi},
\newblock {\it Global maximizers for the sphere adjoint Fourier restriction inequality},
\newblock J. Funct. Anal. {\bf 268} (2015), 690--702.


\bibitem[GN22]{GN22} 
\textsc{F. Gon\c{c}alves and G. Negro},
\newblock {\it  Local maximizers of adjoint Fourier restriction estimates for the cone, paraboloid and sphere},
\newblock  Anal. PDE {\bf 15(4)}: 1097-1130 (2022).

\bibitem[GR14]{GR14} 
\textsc{I. S. Gradshteyn and I. M. Ryzhik},
\newblock Table of Integrals, Series, and Products, 
\newblock (translated from Russian), seventh edition, Elsevier/Academic Press, Amsterdam, 2007.


\bibitem[Kr14]{Kr14}
\textsc{I. Krasikov}, 
\newblock {\it Approximations for the Bessel and Airy functions with an explicit error term,}
\newblock LMS J. Comput. Math. {\bf 17} (2014), no.~1, 209--225. 

\bibitem[{La00}]{La00}
\textsc{L. J. Landau}, 
\newblock {\it Bessel functions: monotonicity and bounds},
\newblock J. London Math. Soc. (2) {\bf 61} (2000), no. 1, 197--215.


\bibitem[NOeST23]{NOeST23}
\textsc{G. Negro, D. Oliveira e Silva, and C. Thiele},
\newblock{\it When does $e^\tau$ maximize Fourier extension for a conic section?},
  Harmonic analysis and convexity, 391--426. Adv. Anal. Geom., 9.


\bibitem[OeSTZK19]{oliveira2019band}
\textsc{D. Oliveira e Silva, C. Thiele, and P. Zorin-Kranich},
\newblock {\it Band-limited maximizers for a {F}ourier extension inequality on the circle,}
\newblock  Exp. Math. {\bf 31} (2022), no. 1, 192--198. 

\bibitem[OeSQ21]{OSQ}
\textsc{D. Oliveira e Silva and R. Quilodr\'an},
\newblock{\it Global mazimizers for adjoint {F}ourier restriction inequalities on low dimensional Spheres.}
\newblock J. Funct. Anal. 280 (2021), no. 7, Paper no. 108825.

\bibitem[{Sh16}]{Sh16}
\textsc{S. Shao},
\newblock {\it On existence of extremizers for the Tomas--Stein inequality for $\mathbb{S}^1$},
\newblock J. Funct. Anal. {\bf 270} (2016), 3996--4038.


\bibitem[{St00}]{St00}
\textsc{K. Stempak},
\newblock {\it A weighted uniform $L^p$-estimate of Bessel functions: a note on a paper of K. Guo},
\newblock Proc. Amer. Math. Soc. {\bf 128} (2000), no. 10, 2943--2945.


\bibitem[Ve88]{Ve88}
\textsc{L. Vega},
\newblock{\it El multiplicador de Schr\"odinger: la funci\'{o}n maximal y los operadores de restricci\'on},
\newblock Ph.D. thesis, Universidad Aut\'onoma de Madrid, 1988.

\bibitem[{Ve92}]{Ve92}
\textsc{L. Vega}, 
\newblock {\it Restriction theorems and the Schr\"{o}dinger multiplier on the torus}, 
\newblock Partial differential equations with minimal smoothness and applications (Chicago, IL, 1990), 199--211, IMA Vol. Math. Appl., 42, Springer, New York, 1992.

\bibitem[{Wa66}]{Wa66}
\textsc{G. N. Watson}, 
\newblock { A Treatise on the Theory of Bessel Functions,}
\newblock Second Edition. Cambridge University Press, Cambridge, 1966.

\end{thebibliography}
\end{document}